%% file: proxtone.tex
\documentclass{article}

\usepackage{algorithm, algpseudocode}
\usepackage{amsfonts, amsmath, amsopn, amsthm, epsfig}
  \numberwithin{equation}{section}
\usepackage{graphicx, epstopdf}
\usepackage{hyperref}
\usepackage{subfigure}
\usepackage{appendix}

\usepackage{calc}
\usepackage[top=2.5cm,bottom=2.3cm,inner=2.3cm,outer=2.3cm,bindingoffset=0.5cm,footskip=0.55cm]{geometry}
\newtheorem{theorem}{Theorem}[section]
\newtheorem{corollary}[theorem]{Corollary}

\theoremstyle{remark}

\newcommand{\E}{\mathbb{E}}

\newenvironment{lemma*}[2][Lemma]{\par\bgroup{\bfseries #1\ #2. }\it\ignorespaces}{\egroup}

\title{Proximal Stochastic Newton-type Gradient Descent Methods for Minimizing Regularized Finite Sums}

\author{Ziqiang Shi\footnotemark[1] \footnotemark[2]
}

\date{August, 2014}

\input{macros.tex}

\begin{document}

\maketitle

\renewcommand{\thefootnote}{\fnsymbol{footnote}}

\footnotetext[1]{Fujitsu Research \& Development Center, Beijing, China.}
\footnotetext[2]{shiziqiang@cn.fujitsu.com; shiziqiang7@gmail.com.}

\renewcommand{\thefootnote}{\arabic{footnote}}

\begin{abstract}
In this work, we generalized and unified recent two completely different works of Jascha~\cite{sohl2014fast} and Lee~\cite{lee2012proximal} respectively into one by proposing the \textbf{prox}imal s\textbf{to}chastic \textbf{N}ewton-type gradient (PROXTONE) method for optimizing the sums of two convex functions: one is the average of a huge number of smooth convex functions, and the other is a non-smooth convex function. While a set of recently proposed proximal stochastic gradient methods, include MISO, Prox-SDCA, Prox-SVRG, and SAG, converge at linear rates, the PROXTONE incorporates second order information to obtain stronger convergence results, that it achieves a linear convergence rate not only in the value of the objective function, but also in the \emph{solution}. The proof is simple and intuitive, and the results and technique can be served as a initiate for the research on the proximal stochastic methods that employ second order information.
\end{abstract}

\section{Introduction and problem statement}
\label{sec:introduction}

In this work, we consider the problems of the following form:
\begin{align}
  \minimize_{x \in \reals^p} \,f(x) := \frac{1}{n}\sum_{i=1}^n g_i (x) + h(x),
  \label{eq:composite-form-online}
\end{align}
where $g_i$ is a smooth convex loss function associated with a sample in a training set, and
$h$ is a non-smooth convex penalty function or
regularizer. Let $g(x)=\frac{1}{n}\sum_{i=0}^n g_i (x)$. We assume the optimal value $f^\star$ is
attained at some optimal solution $x^\star$, not necessarily unique. Problems of this form often arise
in machine learning, such as the least-squares regression, the Lasso, the elastic net, and the logistic regression.

For optimizing~\eqref{eq:composite-form-online}, the standard and popular~\emph{proximal full gradient method} (Prox-FG) uses iterations of the form
\begin{equation}
\label{eqn:prox-fg}
    x^{k+1} = \argmin_{x\in\reals^p} \left\{ \nabla g(x_{k})^T x
    + \frac{1}{2\alpha_k}\|x-x_{k-1}\|^2 + h(x) \right\} ,
\end{equation}
where $\alpha_k$ is the step size at the $k$-th iteration. Under standard
assumptions the sub-optimality achieved on iteration $k$ of the Prox-FG method with a constant step size
is given by
\[
\mathbb{E}[f(x^k)] - f(x^\ast) = O(\frac{1}{k}).
\]
When $f$ is strongly-convex, the error satisfies~\cite{xiao2014proximal}
\[
\mathbb{E}[f(x^k)] - f(x^\ast) = O(\bigl(\frac{L-\mu_g}{L+\mu_h}\bigr)^k),
\]
where $L$ is the Lipschitz constant of $f(x)$, $\mu_g$, and $\mu_h$ are the convexity parameters of $g(x)$ and $h(x)$ respectively. These notations will be detailed in Section~\ref{sec:notations}. This results in a linear convergence rate, which is also known as a geometric or exponential rate
because the error is cut by a fixed fraction on each iteration.

Unfortunately, the Prox-FG and methods can be
unappealing when $n$ is large because its iteration cost scales linearly in $n$.
When the number of components $n$ is very large, each iteration of~\eqref{eqn:prox-fg} can be very
expensive since it requires computing the gradients for all the n component functions $g_i$, and also their average.

The main appeal of proximal stochastic gradient (Prox-SG) methods is that they have an iteration cost which is independent of $n$,
making them suited for modern problems where $n$ may be very large. The basic Prox-SG method for
optimizing~\eqref{eq:composite-form-online}, uses iterations of the form
\begin{equation}\label{eqn:prox-sg}
x_k = \prox_{\alpha_k h}\bigl(x_{k-1} - \alpha_k \nabla g_{i_k}(x_{k-1}) \bigr) ,
\end{equation}
where at each iteration an index $i_k$ is sampled uniformly from the set $\{1, ..., n\}$. The randomly
chosen gradient $\nabla g_{i_k}(x_{k-1})$ yields an unbiased estimate of the true gradient $\nabla g(x_{k-1})$ and one can show
under standard assumptions that,  for a suitably chosen decreasing
step-size sequence $\{\alpha_k\}$, the Prox-SG iterations have an expected sub-optimality for convex objectives of~\cite{bertsekas2011incremental}
\[
\mathbb{E}[f(x^k)] - f(x^\ast) = O(\frac{1}{\sqrt{k}})£¬
\]
and an expected sub-optimality for strongly-convex objectives of
\[
\mathbb{E}[f(x^k)] - f(x^\ast) = O(\frac{1}{k}).
\]
In these rates, the expectations are taken with respect to the selection of the $i_k$ variables.

There is another group of methods, which converges much faster, but need more memory and computation to obtain the second order information about the objective function. These methods are always limited to small-to-medium scale problems that require a high degree of precision. For optimizing~\eqref{eq:composite-form-online}, \emph{proximal Newton-type} methods~\cite{lee2012proximal} that incorporate second order information use iterations of the form $x^{k+1} \leftarrow x^k + \Delta x^k$, here $\Delta x^k$ is obtained by
\begin{equation}
\label{eqn:prox-n}
    \Delta x^k = \argmin_{d\in\reals^p}\, \nabla g(x^k)^Td + \frac{1}{2}d^TH_kd + h(x^k+d),
\end{equation}
where $H_k$ denotes an approximation to $\nabla^2 g(x_k)$.
According to the strategies for choosing $H_k$, we obtain different method, such as \emph{proximal Newton method} (Prox-N) when we choose $H_k$ to be $\nabla^2 g(x^k)$; \emph{proximal quasi-Newton method} (Prox-QN) when we build an approximation to $\nabla^2 g(x_k)$ using
changes measured in $\nabla g$ according to a quasi-Newton strategy~\cite{lee2012proximal}. Indeed if we compared~\eqref{eqn:prox-n} with~\eqref{eqn:prox-fg}, it can be seen Prox-N is the Prox-FG with scaled proximal mappings.

Based on the background above, now we can describe our approaches and findings. The primary contribution of this work is the proposal and analysis of a new algorithm that we call the proximal stochastic Newton-type gradient (PROXTONE, pronounced /prok stone/) method, a stochastic variant of the Prox-N method.
The PROXTONE method has the low iteration cost as that of Prox-SG methods, but achieves
the convergence rates stated above for the Prox-FG method. The PROXTONE iterations take the form $x^{k+1} \leftarrow x^k + t_k\Delta x^k$, where $\Delta x^k$ is obtained by
\begin{equation}
\label{eqn:proxtone}
    \Delta x^k \leftarrow \arg\min_{d} d^T(\nabla_k + H_kx^k)+\frac{1}{2}d^T H_k d+h(x^k+d),
\end{equation}
here $\nabla_{k} = \frac{1}{n}\sum_{i=1}^n \nabla_{k}^i$, $H_{k} =  \frac{1}{n}\sum_{i=1}^n H_{k}^i$, and at each iteration a random index $j$ and corresponding $H_{k+1}^j$ is selected, then we set
\[
\nabla_{k+1}^i = \begin{cases}
 \nabla g_i(x^{k+1})-H_{k+1}^i x^{k+1} & \textrm{if $i = j$,}\\
\nabla_{k+1}^i & \textrm{otherwise.}
\end{cases}
\]
and $H_{k+1}^i\leftarrow H_{k}^i$ ($i\neq j$).

That is, like the Prox-FG and Prox-N method, the step incorporates a gradient with respect to each function. But,
like the Prox-SG method, each iteration only computes the gradient with respect to a single example and
the cost of the iterations is independent of $n$. Despite the low cost of the PROXTONE iterations, we show
in this paper that the PROXTONE iterations have a linear convergence rate for strongly-convex objectives, like the Prox-FG method.
That is, by having access to $j$ and by keeping a memory of the approximation for the Hessian matrix computed
for the objective funtion, this iteration achieves a faster convergence rate than is possible for standard Prox-SG
methods.

There are a large variety of approaches available to accelerate the convergence of Prox-SG methods, and
a full review of this immense literature would be outside the scope of this work. Several recent work considered various special cases
of~\eqref{eq:composite-form-online},
and developed algorithms that enjoy the linear convergence rate, such as Prox-SDCA~\cite{shalev2012proximal}, MISO~\cite{mairal2013optimization}, SAG~\cite{schmidt2013minimizing}, Prox-SVRG~\cite{xiao2014proximal}, SFO~\cite{sohl2014fast}, and Prox-N~\cite{lee2012proximal}. All these methods converge with an exponential rate in the value of the objective function, except that the Prox-N achieves superlinear rates of convergence in \emph{solution}, however it is a batch mode method. Shalev-Shwartz and Zhang~\cite{shalev2013stochastic,shalev2012proximal}'s Prox-SDCA considered the case where the component functions have the form
$g_i(x)=\phi_i(a_i^T x)$ and the Fenchel conjugate functions of~$\phi_i$ and~$h$ can be computed efficiently. Schimidt et al.~\cite{schmidt2013minimizing}'s SAG and Jascha et al.~\cite{sohl2014fast}'s SFO considered the case where $h(x)\equiv 0$.


Our PROXTONE is a extension of the SFO and Prox-N to a proximal stochastic Newton-type method for solving the more \emph{general} ( compared to Prox-SDCA, SAG and SFO) class of problems defined in~\eqref{eq:composite-form-online}. PROXTONE makes connections between two completely different approaches. It achieves a linear convergence rate not only in the value of the objective function, but also in the \emph{solution}. We now outline the rest of the study. Section~\ref{sec:PROXTONE} presents the main algorithm and gives a equivalent form in order for the ease of analysis.
Section~\ref{sec:Analysis} states the assumptions underlying our analysis and gives the main results; we first give a linear convergence rate in \emph{function value} (weak convergence) that applies for any problem, and then give a strong linear convergence rate in \emph{solution}, however with some additional conditions. Finally we conclude in Section~\ref{sec:Conclusions}.


\subsection{Notations and Assumptions}
\label{sec:notations}

Before proceeding, we introduce the notations and some useful lemmas formally first. In this work, we most adopt the nomenclature used by Nesterov~\cite{Nesterov04book}. The functions encountered in this work are all convex if there are no other statements.

In this paper, we assume the function $h(x)$ is lower semi-continuous and convex,
and its effective domain, $\dom(h):=\{x\in\reals^p\,|\,h(x)<+\infty\}$,
is closed. Each $g_i(x)$, for $i=1,\ldots,n$, is differentiable
on an open set that contains~$\dom(h)$, and their gradients are
Lipschitz continuous.
That is, there exist~$L_i>0$ such that for all $x, y\in\dom(h)$,
\begin{equation}\label{eqn:smooth-i}
    \|\nabla g_i(x) - \nabla g_i(y)\| \leq L_i \|x-y\|.
\end{equation}
This is a fairly weak assumption on the $g_i$ functions, and in cases where the $g_i$ are twice-differentiable
it is equivalent to saying that the eigenvalues of the Hessians of each $g_i$ are bounded above by $L_i$.
Then from the Lemma 1.2.3 and its proof in Nesterov's book, for $i=1,\ldots,n$, we have
\begin{equation}\label{eqn:smooth-i-1}
    | g_i(x) -  g_i(y) - \nabla g_i(y)^T(x-y)| \leq \frac{L_i}{2} \|x-y\|^2.
\end{equation}
The above assumption of~\eqref{eqn:smooth-i} implies that the gradient
of the average function $g(x)$ is also Lipschitz continuous, i.e.,
there is an $L>0$ such that for all $x, y\in\dom(h)$,
\[
    \|\nabla g(x) - \nabla g(y)\| \leq L \|x-y\|.
\]
Moreover, we have $L\leq (1/n)\sum_{i=1}^n L_i$.

A function $f(x)$ is called $\mu$-strongly convex, if
there exist $\mu\ge 0$ such that for all $x\in\dom(f)$ and $y\in\reals^p$,
\begin{equation}\label{eqn:strong-convex}
    f(y) \geq f(x) + \xi^T(y-x) + \frac{\mu}{2} \|y-x\|^2,
    \quad \forall\, \xi\in\partial f(x).
\end{equation}
The \emph{convexity parameter} of a function is the largest~$\mu$
such that the above condition holds. If $\mu = 0$, it is identical to the definition of a convex function.
The strong convexity of $f(x)$ in~\eqref{eq:composite-form-online} may come from either $g(x)$ or $h(x)$ or both.
More precisely, let $g(x)$ and  $h(x)$ have convexity parameters $\mu_g$
and $\mu_h$ respectively, then $\mu \geq \mu_g + \mu_h$.
From Lemma B.5 in~\cite{mairal2013optimization} and~\eqref{eqn:strong-convex}, we have
\begin{equation}\label{eqn:(Second-Order-Growth}
    f(y) \geq f(x^*)  + \frac{\mu}{2} \|y-x^*\|^2.
\end{equation}

\section{The PROXTONE method}
\label{sec:PROXTONE}

We summarize the PROXTONE method of~\eqref{eqn:proxtone} in Algorithm~\ref{alg:General-proxtone}. It can be easily checked that if $n=1$, then it becomes the determined proximal Newton-type methods proposed by Lee and Sun et al.~\cite{lee2012proximal} for minimizing composite functions:
\begin{align}
  \minimize_{x \in \reals^p} \,f(x) := g(x) + h(x)
\end{align}
by~\eqref{eqn:prox-n},
thus PROXTONE is indeed a generalization of Prox-NG.

\begin{algorithm}[H]
\caption{PROXTONE: A generic PROXimal sTOchastic NEwton-type gradient descent method}
\label{alg:General-proxtone}

\textbf{Input}: start point $x^0 \in$ dom $f$; for $i\in\{1,2,..,n\}$, let $H_{-1}^i=H_0^i$ be a positive definite approximation to the Hessian of $g_i(x)$ at $x^0$, $\nabla_{-1}^i=\nabla_0^i=\nabla g_i(x^0) - H_0^ix^0$; and $\nabla_0=\frac{1}{n}\sum_{i=1}^n \nabla_0^i$, $H_0=\frac{1}{n}\sum_{i=0}^n H_0^i$.

1: \textbf{repeat}

2: Solve the subproblem for a search direction:$\triangle x^k \leftarrow \arg\min_{d} d^T(\nabla_k + H_kx^k)+\frac{1}{2}d^T H_k d+h(x^k+d)$.


3: Update: $x^{k+1}=x^{k}+\triangle x^k$.

4: Sample $j$ from $\{1,2,..,n\}$, use the $\nabla g_j(x^{k+1})$ and $H_{k+1}^j$, which is a positive definite approximation to the Hessian of $g_j(x)$ at $x^{k+1}$, to update the $\nabla_{k+1}^i$ ($i\in\{1,2,..,n\}$):  $\nabla_{k+1}^j\leftarrow \nabla g_j(x^{k+1})-H_{k+1}^j x^{k+1}$, while leaving all other $\nabla_{k+1}^i$ and $H_{k+1}^i$ unchanged: $\nabla_{k+1}^i\leftarrow \nabla_{k}^i$ and $H_{k+1}^i\leftarrow H_{k}^i$ ($i\neq j$) ; and finally obtain $\nabla_{k+1}$ and $H_{k+1}$ by $\nabla_{k+1} \leftarrow \frac{1}{n}\sum_{i=1}^n \nabla_{k+1}^i$, $H_{k+1} \leftarrow  \frac{1}{n}\sum_{i=1}^n H_{k+1}^i$.

5: \textbf{until} stopping conditions are satisfied.

\textbf{Output}: $x^k$.
\end{algorithm}

 It is also a generalization of recent work by Jascha~\cite{sohl2014fast}, whose SFO is the special case of our PROXTONE with $h(x)\equiv 0$. Our algorithm in Jascha's style is summarized in Algorithm~\ref{alg:General-prostone-analysis} which is equivalent to the original PROXTONE. To see the equivalence, keep in mind that $G^k(x)$ is a quadratic function, we only need to check the following equations:
 \[
    \nabla^2G^{k}(x) = \frac{1}{n}\sum_{i=1}^n H_k^i\ \ \text{and}\ \ \nabla G^{k}(x) = \frac{1}{n}\sum_{i=1}^n \nabla g_i(x)+\frac{1}{n}\sum_{i=1}^n (x-x^k)^TH_k^i,
\]
and
\[
\nabla_k + H_kx^k = \frac{1}{n}\sum_{i=1}^n[\nabla g_i(x^{\theta_{i,k}})+(x^k-x^{\theta_{i,k-1}})^TH_{\theta_{i,k}}^i].
\]

In following analysis of Section~\ref{sec:Analysis}, we will not distinguish these two forms from each other.

\begin{algorithm}[H]
\caption{PROXTONE in a form that is easy to analyze}
\label{alg:General-prostone-analysis}

\textbf{Input}: start point $x^0 \in$ dom $f$; for $i\in\{1,2,..,n\}$, let $g_i^0(x)=g_i(x^0)+(x-x^0)^T\nabla g_i(x^0)+\frac{1}{2}(x-x^0)^TH_0^i(x-x^0)$, where the notation $H_0^i$ ($i\in\{1,2,..,n\}$) are totally the same as they in Algorithm~\ref{alg:General-proxtone}; and $G^0(x)=\frac{1}{n}\sum_{i=1}^n g_i^0(x)$.

1: \textbf{repeat}

2: Solve the subproblem for new approximation of the solution: $x^{k+1} \leftarrow \arg\min_{x} \bigl[ G^k(x) + h(x) \bigr] $.

3: Sample $j$ from $\{1,2,..,n\}$, and update the surrogate functions:
\begin{align}
g_j^{k+1}(x)=g_j(x^{k+1})+(x-x^{k+1})^T\nabla g_j(x^{k+1})+\frac{1}{2}(x-x^{k+1})^TH_{k+1}^i(x-x^{k+1}),
  \label{eq:subfunction-surrogate-update}
\end{align}
while leaving all other $g_i^{k+1}(x)$ unchanged: $g_i^{k+1}(x)\leftarrow g_i^{k}(x)$ ($i\neq j$); and $G^{k+1}(x)=\frac{1}{n}\sum_{i=1}^n g_i^{k+1}(x)$.

4: \textbf{until} stopping conditions are satisfied.

\textbf{Output}: $x^k$.
\end{algorithm}

\section{Convergence Analysis}
\label{sec:Analysis}

Under the standard assumptions, we now state our convergence result.

\begin{theorem}\label{thm:prox-stone}
Suppose $\nabla g_i(x)$ is Lipschitz continuous with constant $L_i > 0$ for $i=1,...,n$, and $\L_i I \preceq mI\preceq H_{k}^i \preceq MI$ for all $i=1,...,n$ and $k\geq1$, $h(x)$ is strongly convex with $\mu_h\geq 0$, then the PROXTONE iterations satisfy for $k\geq 1$:
%
\begin{align}
   \mathbb{E}[f(x^k)] - f^* \le \frac{M+L_{max}}{2} [\frac{1}{n} \frac{M+L_{max}}{2\mu_h+m}+(1-\frac{1}{n})]^k\|x^*-x^0\|^2.
\label{eq:proxtone-value-linear}
\end{align}
\end{theorem}

The ideas of the proof is near identical to that of MISO by Mairal~\cite{mairal2013optimization} and for completeness we give a simple version in the appendix.

We have the following remarks regarding the above result:
\begin{itemize}
\item
In order to satisfy $\mathbb{E}[f(x^k)] - f^* \leq \epsilon$,
the number of iterations~$k$ needs to satisfy
\[
    k\geq (\log\rho)^{-1} \log \bigl[\frac{2\epsilon}{(M+L_{max})\|x^*-x^0\|^2}\bigr],
\]
where $\rho=\frac{1}{n} \frac{M+L_{max}}{2\mu_h+m}+(1-\frac{1}{n})$.
\item
Inequality~\eqref{eq:proxtone-value-linear} gives us a reliable stopping criterion for the PROXTONE method.
\end{itemize}

At this moment, we see that the expected quality of the output of PROXTONE is good. However, in practice we are not going to run this method many times on the same problem. What is the probability that our single run can give us also a good result.
Since $f(x^k) - f^*  \geq 0$, Markov's
inequality and Theorem~\ref{thm:prox-stone} imply
that for any $\epsilon>0$,
\[
\Prob\Bigl( f(x^k) - f^* \geq \epsilon\Bigr)
~\leq~ \frac{\E [f(x^k) - f^*]}{\epsilon}
~\leq~ \frac{(M+L_{max})\rho^k\|x^*-x^0\|^2 }{2\epsilon}.
\]
Thus we have the following high-probability bound.

\begin{corollary}\label{cor:high-prob-1}
Suppose the assumptions in Theorem~\ref{thm:prox-stone} hold.
Then for any $\epsilon>0$ and $\delta\in(0,1)$, we have
\[
  \Prob\bigl(f(x^k) - f(x^\star) \leq \epsilon \bigr) \geq 1-\delta
\]
provided that the number of iterations~$k$ satisfies
\[
k \geq \log\left(\frac{(M+L_{max})\|x^*-x^0\|^2}{2\delta\epsilon}\right)
\bigg/\log\left(\frac{1}{\rho}\right).
\]
\end{corollary}

Based on Theorem~\ref{thm:prox-stone} and its proof, we give a deeper and stronger result that the PROXTONE achieves a linear convergence rate in the solution.

\begin{theorem}\label{thm:prox-stone-linear}
Suppose $\nabla g_i(x)$ and $\nabla^2g_i$ are Lipschitz continuous with constant $L_i > 0$ and $K_i > 0$ respectively for $i=1,...,n$, $h(x)$ is strongly convex with $\mu_h\geq 0$. If $H_{\theta_{i,k}}^i=\nabla^2g_i(x^{\theta_{i,k}})$ and $\L_i I \preceq mI\preceq H_{k}^i \preceq MI$, then PROXTONE converges exponentially to $x^\star$ in expectation:
\[
  \mathbb{E}[ \norm{x^{k+1} - x^\star}] \le (\frac{K_{avg}+2L_{max}}{m} \frac{M+L_{max}}{2\mu_h+m}+\frac{2L_{max}}{m})[\frac{1}{n} \frac{M+L_{max}}{2\mu_h+m}+(1-\frac{1}{n})]^{k-1} \|x^*-x^0\|^2.
\]
\end{theorem}


In order to satisfy $ \mathbb{E}[ \norm{x^{k+1} - x^\star}] \leq \epsilon$,
the number of iterations~$k$ needs to satisfy
\[
    k\geq (\log\rho)^{-1} \log \bigl[\frac{\epsilon}{C\|x^*-x^0\|^2}\bigr],
\]
where $\rho$ is as before and $C=\frac{K_{avg}+2L_{max}}{m} \frac{M+L_{max}}{2\mu_h+m}+\frac{2L_{max}}{m}$.

Due to the Markov's
inequality, Theorem~\ref{thm:prox-stone-linear} implies the following result.

\begin{corollary}\label{cor:high-prob-2}
Suppose the assumptions in Theorem~\ref{thm:prox-stone-linear} hold.
Then for any $\epsilon>0$ and $\delta\in(0,1)$, we have
\[
  \Prob\bigl(\norm{x^{k+1} - x^\star} \geq \epsilon \bigr) \geq 1-\delta
\]
provided that the number of iterations~$k$ satisfies
\[
k \geq \log\left(\frac{((K_{avg}+2L_{max}) (M+L_{max})+2L_{max}(2\mu_h+m)) \|x^*-x^0\|^2}{m(2\mu_h+m)\delta\epsilon}\right)
\bigg/\log\left(\frac{1}{\rho}\right).
\]
\end{corollary}

\section{Conclusions}
\label{sec:Conclusions}

This paper introduces a proximal stochastic method called PROXTONE for minimizing regularized finite sums. For smooth and strongly convex problems, we show that PROXTONE not only enjoys the same linear rates as those of MISO, SAG, Prox-SVRG and Prox-SDCA, but also prove that the \emph{solution} of this method converges in exponential rate too. There are some directions that the current study can be extended. In this paper, we have focused on the theory of PROXTONE; it would be meaningful to also do the numerical evaluation and implementation details~\cite{sohl2014fast}. Second, combine with randomized block coordinate method~\cite{nesterov2012efficiency} for minimizing regularized convex functions with a huge number of varialbes/coordinates. Moreover, due to the trends and needs of big data, we are designing distributed/parallel PROXTONE for real life applications. In a broader context, we believe that the current paper could serve as a basis for examining the method on the proximal stochastic methods that employ second order information.


\appendix
\section*{Appendix}

In this Appendix, we give the proofs of the two propositions.

\section{Proof of Theorem~\ref{thm:prox-stone}}

Since in each iteration of the PROXTONE, we obtain a quadratic function $g_i^{k}(x)$ with random parameters to approximate each $g_i(x)$:
\begin{align}
g_i^{k}(x)=g_i(x^{\theta_{i,k}})+(x-x^{\theta_{i,k}})^T\nabla g_i(x^{\theta_{i,k}})+\frac{1}{2}(x-x^{\theta_{i,k}})^TH_{\theta_{i,k}}^i(x-x^{\theta_{i,k}}),
  \label{eq:subfunction-surrogate-assum}
\end{align}
where $\theta_{i,k}$ is a random variable which have the following conditional probability distribution in each iteration:
\begin{align}
\mathbb{P}(\theta_{i,k}=k|j)=\frac{1}{n} \quad \text{and} \quad \mathbb{P}(\theta_{i,k}=\theta_{i,k-1}|j)=1-\frac{1}{n},
  \label{eq:random-parameter-dist}
\end{align}
that yields
\begin{align}
\label{eq:x-update-relation}
\mathbb{E}[\|x^*-x^{\theta_{i,k}}\|^2]=\frac{1}{n}\mathbb{E}[\|x^*-x^k\|^2]+(1-\frac{1}{n})\mathbb{E}[\|x^*-x^{\theta_{i,k-1}}\|^2].
\end{align}

Since $0\preceq H^i_{\theta_{i,k}} \preceq MI$ and $\nabla^2 g_i^{k}(x)=H_{\theta_{i,k}}^i$, by Theorem~2.1.6 of~\cite{Nesterov04book} and the assumption, $\nabla g_i^{k}(x)$ and $\nabla g_i(x)$ are Lipschitz continuous with constant $M$ and $L_i$ respectively, and further $\nabla g_i^{k}(x)-\nabla g_i(x)$ is Lipschitz continuous with constant $M+L_i$ for $i=1,\ldots,n$. This together with~\eqref{eqn:smooth-i-1} yieds
\[
    | [g_i^{k}(x) - g_i(x)] -  [g_i^{k}(y) - g_i(y)] - \nabla [g_i^{k}(y) - g_i(y)]^T(x-y)| \leq \frac{M+L_i}{2} \|x-y\|^2.
\]
Applying the above inequality with $y=x^{\theta_{i,k}}$, and using the fact that $\nabla [g_i^{k}(x^{\theta_{i,k}})] = \nabla [g_i(x^{\theta_{i,k}})]$ and $g_i^{k}(x^{\theta_{i,k}}) = g_i(x^{\theta_{i,k}})$, we have
\[
    | g_i^{k}(x) - g_i(x)| \leq \frac{M+L_i}{2} \|x-x^{\theta_{i,k}}\|^2.
\]
Summing over $i=1,\ldots,n$ yields
\begin{align}
\label{eq:diff_bound}
    [G^{k}(x)+h(x)] - [g(x)+h(x)] \leq \frac{1}{n}\sum_{i=1}^n\frac{M+L_i}{2} \|x-x^{\theta_{i,k}}\|^2.
\end{align}
Then by the Lipschitz continuity of $\nabla g_i(x)$ and the assumption $L_i I \preceq mI\preceq H_{k}^i$, we have
\begin{eqnarray*}
   g_i(x) &\leq&  g_i(x^{\theta_{i,k}}) + \nabla g_i(x^{\theta_{i,k}})^T(x-x^{\theta_{i,k}})| + \frac{L_i}{2} \|x-x^{\theta_{i,k}}\|^2 \\
   &\leq& g_i(x^{\theta_{i,k}})+(x-x^{\theta_{i,k}})^T\nabla g_i(x^{\theta_{i,k}})+\frac{1}{2}(x-x^{\theta_{i,k}})^TH_{\theta_{i,k}}^i(x-x^{\theta_{i,k}})=g_i^{k}(x),
\end{eqnarray*}
and thus, by summing over $i$ yields $ g(x) \leq G^{k}(x)$, and further by the optimality of $x^{k+1}$, we have
\begin{align}
\label{eq:fx_bound}
    f(x^{k+1})\leq G^{k}(x^{k+1})+h(x^{k+1})\leq G^{k}(x)+h(x) \leq  f(x) + \frac{1}{n}\sum_{i=1}^n\frac{M+L_i}{2} \|x-x^{\theta_{i,k}}\|^2
\end{align}
Since $mI\preceq H_{\theta_{i,k}} $ and $\nabla^2 g_i^{k}(x)=H_{\theta_{i,k}}$, by Theorem~2.1.11 of~\cite{Nesterov04book}, $g_i^{k}(x)$ is $m$-strongly convex. Since $G^{k}(x)$ is the average of $g_i^{k}(x)$, thus $G^{k}(x)+h(x)$ is ($m+\mu_h$)-strongly convex, we have
\begin{eqnarray*}
f(x^{k+1})+\frac{m+\mu_h}{2}\|x-x^{k+1}\|^2
&\leq& G^{k}(x^{k+1})+h(x^{k+1})+\frac{m+\mu_h}{2}\|x-x^{k+1}\|^2 \\
&\leq& G^{k}(x)+h(x) \\
&=& f(x)+[G^{k}(x)+h(x)-f(x)]\\
&\leq& f(x)+\frac{1}{n}\sum_{i=1}^n\frac{M+L_i}{2} \|x-x^{\theta_{i,k}}\|^2.
\end{eqnarray*}
By taking the expectation of both sides and let $x=x^*$ yields
\[
   \mathbb{E}[f(x^{k+1})] - f^* \leq \mathbb{E}[\frac{1}{n}\sum_{i=1}^n\frac{M+L_i}{2} \|x^*-x^{\theta_{i,k}}\|^2]-\mathbb{E}[\frac{m+\mu_h}{2}\|x^*-x^{k+1}\|^2].
\]

%
We have
\[
    \frac{\mu_h}{2} \|x^{k+1}-x^*\|^2  \leq \mathbb{E}[f(x^{k+1})] - f^* \leq  \mathbb{E}[\frac{1}{n}\sum_{i=1}^n\frac{M+L_{max}}{2} \|x-x^{\theta_{i,k}}\|^2]-\mathbb{E}[\frac{m+\mu_h}{2}\|x-x^{k+1}\|^2].
\]
thus
\begin{align}
\label{eq:x_diff_bound}
     \|x^{k+1}-x^*\|^2   \leq  \frac{M+L_{max}}{2\mu_h+m}\mathbb{E}[\frac{1}{n}\sum_{i=1}^n\|x^*-x^{\theta_{i,k}}\|^2].
\end{align}
then we have
\begin{eqnarray*}
  \mathbb{E}[\frac{1}{n}\sum_{i=1}^n \|x^*-x^{\theta_{i,k}}\|^2]
   &=&  \frac{1}{n}\|x^k-x^*\|^2+(1-\frac{1}{n})\mathbb{E}[\frac{1}{n}\sum_{i=1}^n \|x^*-x^{\theta_{i,k-1}}\|^2]  \\
  &\leq&  \frac{1}{n}\|x^k-x^*\|^2+(1-\frac{1}{n})\mathbb{E}[\frac{1}{n}\sum_{i=1}^n \|x^*-x^{\theta_{i,k-1}}\|^2] \\
   &\leq&  [\frac{1}{n} \frac{M+L_{max}}{2\mu_h+m}+(1-\frac{1}{n})]\mathbb{E}[\frac{1}{n}\sum_{i=1}^n \|x^*-x^{\theta_{i,k-1}}\|^2] \\
    &\leq&  [\frac{1}{n} \frac{M+L_{max}}{2\mu_h+m}+(1-\frac{1}{n})]^k\mathbb{E}[\frac{1}{n}\sum_{i=1}^n \|x^*-x^{\theta_{i,0}}\|^2] \\
      &\leq&  [\frac{1}{n} \frac{M+L_{max}}{2\mu_h+m}+(1-\frac{1}{n})]^k \|x^*-x^0\|^2 .
\end{eqnarray*}

Thus we have  $\mathbb{E}[f(x^{k+1})] - f^* \leq \frac{M+L_{max}}{2} [\frac{1}{n} \frac{M+L_{max}}{2\mu_h+m}+(1-\frac{1}{n})]^k \|x^*-x^0\|^2$.

\section{Proof of Theorem~\ref{thm:prox-stone-linear}}

We first examine the relations between the search directions of Prox-N and PROXTONE.

By \eqref{eqn:prox-n}, \eqref{eqn:proxtone} and Fermat's rule, $\Delta x^k_{Prox-N}$ and $\Delta x^k$ are also the solutions to
\begin{align*}
\Delta x^k_{Prox-N} &= \argmin_{d\in\reals^p}\, \nabla g(x^k)^Td + (\Delta x^k_{Prox-N})^TH_kd + h(x^k+d), \\
\Delta x^k  &= \argmin_{d\in\reals^p}\, (\nabla_k + H_kx^k)^Td + (\Delta x^k)^TH_kd + h(x^k+d).
\end{align*}
Hence $\Delta x^k$ and $\Delta x^k_{Prox-N}$ satisfy
\begin{align*}
&\nabla g(x^k)^T\Delta x^k + (\Delta x^k_{Prox-N})^TH_k\Delta x^k + h(x^k+\Delta x^k) \\
\geq &\nabla g(x^k)^T\Delta x^k_{Prox-N} + (\Delta x^k_{Prox-N})^TH_k\Delta x^k_{Prox-N} + h(x^k+\Delta x^k_{Prox-N})
\end{align*}
and
\begin{align*}
 &(\nabla_k + H_kx^k)^T\Delta x^k_{Prox-N} + (\Delta x^k)^TH_k\Delta x^k_{Prox-N} + h(x^k+\Delta x^k_{Prox-N}) \\
\geq &(\nabla_k + H_kx^k)^T\Delta x^k + (\Delta x^k)^TH_k\Delta x^k + h(x^k+\Delta x^k).
\end{align*}
We sum these two inequalities and rearrange to obtain
\begin{align*}
&(\Delta x^k)^TH_k\Delta x^k - 2(\Delta x^k_{Prox-N})^TH_k\Delta x^k + (\Delta x^k_{Prox-N})^TH_k\Delta x^k_{Prox-N} \\ &\le  (\nabla_k + H_kx^k - \nabla g(x^k))^T(\Delta x^k_{Prox-N}-\Delta x^k ).
\end{align*}

The assumptions $mI\preceq H_{\theta_{i,k}} $ yields that $mI\preceq H_k $, together with we have
\[
m\|\Delta x^k - \Delta x^k_{Prox-N}\|^2 \le \| \frac{1}{n}\sum_{i=1}^n(\nabla g_i(x^{\theta_{i,k}})- \nabla g_i(x^k)-(x^{\theta_{i,k}}-x^k)^TH_{\theta_{i,k}}^i)\|\|(\Delta x^k - \Delta x^k_{Prox-N})\|.
\]
Since we have
\[
\|\Delta x^k - \Delta x^k_{Prox-N}\| \le \frac{K_{max}}{2mn} \sum_{i=1}^n \|x^{\theta_{i,k-1}}-x^k\|^2
\]

Since the Prox-N method converges $q$-quadratically (\cf\ Theorem 3.3 of~\cite{lee2012proximal}),
\begin{align*}
\norm{x^{k+1} - x^\star} &\le \norm{x^k + \Delta x^k_{Prox-N} - x^\star}
                         + \norm{\Delta x^k - \Delta x^k_{Prox-N}} \nonumber \\
&\le \frac{K_{avg}}{m}\norm{x^k-x^\star}^2 + \norm{\Delta x^k
  - \Delta x^k_{Prox-N}},
\end{align*}
where $\Delta x^k_{Prox-N}$ denotes the Prox-N search direction.


Combine by we have almost surely that
\begin{align*}
\norm{x^{k+1} - x^\star} &\le \frac{L_2}{m}\norm{x^k-x^\star}^2 + \frac{L_{max}}{2mn} \sum_{i=1}^n \|x^{\theta_{i,k-1}}-x^k\|^2 \\
&\le \frac{K_{avg}}{m}\norm{x^k-x^\star}^2 + \frac{L_{max}}{mn} \sum_{i=1}^n 2\|x^{\theta_{i,k-1}}-x^*\|^2+\frac{L_{max}}{mn} \sum_{i=1}^n 2\|x^*-x^k\|^2.
\end{align*}
Then by~\eqref{eq:x_diff_bound}, we have
\[
\norm{x^{k+1} - x^\star} \le (\frac{K_{avg}+2L_{max}}{m} \frac{M+L_{max}}{2\mu_h+m}+\frac{2L_{max}}{m})\mathbb{E}[\frac{1}{n}\sum_{i=1}^n\|x^{\theta_{i,k}}-x^*\|^2]
\]
which yieds
\[
\norm{x^{k+1} - x^\star} \le (\frac{K_{avg}+2L_{max}}{m} \frac{M+L_{max}}{2\mu_h+m}+\frac{2L_{max}}{m})[\frac{1}{n} \frac{M+L_{max}}{2\mu_h+m}+(1-\frac{1}{n})]^k \|x^*-x^0\|^2.
\]

\frenchspacing
\bibliographystyle{plain}
\bibliography{proxtone_arxiv}

\end{document}

%% file: macros.tex
\newcommand{\BALD}{\begin{aligned}}
\newcommand{\EALD}{\end{aligned}}
\newcommand{\BALDS}{\begin{aligned*}}
\newcommand{\EALDS}{\end{aligned*}}
\newcommand{\BCAS}{\begin{cases}}
\newcommand{\ECAS}{\end{cases}}
\newcommand{\BEAS}{\begin{eqnarray*}}
\newcommand{\EEAS}{\end{eqnarray*}}
\newcommand{\BEQ}{\begin{equation}}
\newcommand{\EEQ}{\end{equation}}
\newcommand{\BIT}{\begin{itemize}}
\newcommand{\EIT}{\end{itemize}}
\newcommand{\BMAT}{\begin{bmatrix}}
\newcommand{\EMAT}{\end{bmatrix}}
\newcommand{\BNUM}{\begin{enumerate}}
\newcommand{\ENUM}{\end{enumerate}}

\newcommand{\cf}{{\it cf.}}

\newcommand{\BA}{\begin{array}}
\newcommand{\EA}{\end{array}}

\newcommand{\reals}{\mathbf{R}}

\DeclareMathOperator*{\argmin}{\arg\min}

\DeclareMathOperator*{\minimize}{minimize}

\newcommand{\Prob}{\mathop{\mathbf{Pr}}}

\newcommand{\norm}[1]{\left\| #1 \right\|}

\DeclareMathOperator{\dom}{dom}

\DeclareMathOperator{\prox}{prox}